\documentclass[11pt,a4paper,leqno]{article}

\usepackage{amsmath}
\usepackage{amsfonts}
\usepackage{amsthm}
\usepackage{amssymb}

\title{Divisorial contractions in dimension $3$\\
which contract divisors to smooth points}
\author{Masayuki Kawakita}
\date{}
\theoremstyle{definition}
\newtheorem{defn}{Definition}[section]
\newtheorem{constr}[defn]{Construction}

\theoremstyle{plain}
\newtheorem{thm}[defn]{Theorem}
\newtheorem{prop}[defn]{Proposition}
\newtheorem{lem}[defn]{Lemma}

\theoremstyle{remark}
\newtheorem{rem}{Remark}[defn]
\newtheorem{cl}[defn]{Claim}

\newtheoremstyle{citing}
{3pt}{3pt}{\itshape}{}{\bfseries}{.}{.5em}{\thmnote{#3}}

\theoremstyle{citing}
\newtheorem*{varthm}{}

\newcommand{\sO}{{\mathcal O}}
\newcommand{\sQ}{{\mathcal Q}}
\newcommand{\mP}{{\mathbb P}}
\newcommand{\mQ}{{\mathbb Q}}
\newcommand{\mZ}{{\mathbb Z}}
\newcommand{\Spec}{\mathrm{Spec}}
\newcommand{\Proj}{\mathrm{Proj}}
\newcommand{\wdiv}{\mathrm{div}}
\newcommand{\Span}{\mathrm{Span}}
\newcommand{\gm}{{\mathfrak m}}

\numberwithin{equation}{section}

\begin{document}

\maketitle

\begin{abstract}
We deal with a divisorial contraction in dimension $3$ which
contracts its exceptional divisor to a smooth point.
We prove that any such contraction can be obtained by a suitable
weighted blow-up.
\end{abstract}

\setcounter{section}{-1}
\section{Introduction}
Divisorial contractions play a major role
in the minimal model program (\cite{KMM87}).
Now that we know this program works in dimension $3$ (\cite{M88}),
it is desirable to describe them explicitly in dimension $3$.
Moreover also in view of the Sarkisov program (\cite{Co95}) and
its applications (for example \cite{CPR99}),
we can recognize the importance of such description
since Sarkisov links of types I and II in this program start from
the converse of divisorial contractions.

Now we concentrate on divisorial contractions in dimension $3$.
Let $f \colon (Y \supset E) \to (X \ni P)$ be such a contraction.
There are two ways to deal with $f$, that is to say, one starting
from $Y$, and the other from $X$.
From the former standpoint,
S.\@ Mori classified them in the case when $Y$ is smooth (\cite{M82}),
and S.\@ Cutkosky extended this result to the case
when $Y$ has only terminal Gorenstein singularities (\cite{Cu88}).
On the other hand, from the latter standpoint, Y.\@ Kawamata showed
that $f$ must be a certain weighted blow-up when $P$ is a terminal quotient
singularity (\cite{K96}), and A.\@ Corti showed that $f$ must be the
blow-up when $P$ is an ordinary double point (\cite[Theorem 3.10]{Co99}).

While it seems that singularities on $Y$ make it hard to tackle the problem
in the former case, the singularity of $P$ may be useful in the latter case
because it gives a special filtration in the tangent space at $P$.
In this paper we treat the case when $P$
is a smooth point and prove the following theorem:

\begin{varthm}[Theorem \ref{mainthm}]
Let $Y$ be a $3$-dimensional $\mQ$-factorial normal
variety with only terminal singularities,
and let $f \colon (Y \supset E) \to (X \ni P)$ be an algebraic germ of
a divisorial contraction which contracts its exceptional divisor
$E$ to a smooth point $P$. 
Then we can take local parameters $x, y, z$ at $P$
and coprime positive integers $a$ and $b$,
such that $f$ is the weighted blow-up of $X$ with its weights
$(x, y, z) = (1, a, b)$.
\end{varthm}

Now we explain our approach to the problem. Y.\@ Kawamata 
adopted the method of comparing discrepancies
of exceptional divisors, and A.\@ Corti applied Shokurov's
connectedness lemma (\cite[Theorem 17.4]{K+92}).
But in the case when $P$ is a smooth point, these methods do not
work well if the center of $E$ on $\mathrm{Bl}_P(X)$ is a point.
Our main tools are the singular Riemann-Roch formula
(\cite[Theorem 10.2]{R87}) on $Y$ and a relative vanishing theorem
(\cite[Theorem 1-2-5]{KMM87}) with respect to $f$.
First with them we derive a rather simple formula
for $\dim_k\sO_X/f_*\sO_Y(-iE)$'s and
an upper-bound of the number of fictitious non-Gorenstein points
of $Y$ (Proposition \ref{eqs}). 
Next using this upper-bound,
we show that the coefficient of $E$ in the pull-back of
a general prime divisor through $P$ is $1$ (Subsection \ref{2E}).
And finally investigating the values of
$\dim_k\sO_X/f_*\sO_Y(-iE)$'s more carefully,
we prove the theorem (Subsection \ref{nE}).

I wish to express my gratitude to
Professor Yujiro Kawamata for his valuable comments and warm encouragement.
He also recommended me to read the papers \cite{CPR99} and \cite{Co99}.
In fact I found the problem treated here as \cite[Conjecture 3.11]{Co99}.

\section{Statement of the theorem}
We work over an algebraically closed field $k$ of characteristic zero.
A variety means an integral separated scheme of finite type over $\Spec \, k$.
We use basic terminologies in \cite[Chapters 1,\,\,2]{K+92}.

Before we state the theorem, we have to define a divisorial contraction.
In this paper it means a morphism which may emerge in
the minimal model program (see \cite{KMM87}).

\begin{defn}\label{divisorial}
Let $f \colon Y \to X$ be a morphism with connected fibers between
normal varieties. We call $f$ a {\itshape divisorial contraction}
if it satisfies the following conditions:
\begin{enumerate}
\item Y is $\mQ$-factorial with only terminal singularities.
\item The exceptional locus of $f$ is a prime divisor.
\item $-K_Y$ is $f$-ample.
\item The relative Picard number of $f$ is $1$.
\end{enumerate}
\end{defn}

Now it is the time when we state the theorem precisely.

\begin{thm}\label{mainthm}
Let $Y$ be a $3$-dimensional $\mQ$-factorial normal
variety with only terminal singularities,
and let $f \colon (Y \supset E) \to (X \ni P)$ be an algebraic germ of
a divisorial contraction which contracts its exceptional divisor
$E$ to a smooth point $P$. 
Then we can take local parameters $x, y, z$ at $P$
and coprime positive integers $a$ and $b$,
such that $f$ is the weighted blow-up of $X$ with its weights
$(x, y, z) = (1, a, b)$.
\end{thm}

\section{Proof of the theorem}
\subsection{Strategy for its proof}
We may assume that $X$ is projective and smooth,
and consider its algebraic germ if necessary. 
First we construct a series of birational morphisms.

\begin{constr}\label{constr}
We construct birational morphisms $g_i \colon  X_i \to X_{i-1}$
between smooth varieties,
integral closed subschemes $Z_i \subset X_i$,
and prime divisors $F_i$ on $X_i$ inductively,
and define positive integers $n, m$, with the following procedure:

\begin{enumerate}
\item Define $X_0$ as $X$ and $Z_0$ as $P$.
\item Let $b_i \colon \mathrm{Bl}_{Z_{i-1}}(X_{i-1}) \to X_{i-1}$ be
the blow-up of $X_{i-1}$ along $Z_{i-1}$, and let ${b'}_{\!i} \colon X_i \to
\mathrm{Bl}_{Z_{i-1}}(X_{i-1})$ be a resolution of
$\mathrm{Bl}_{Z_{i-1}}(X_{i-1})$, that is,
a proper birational morphism from a smooth variety $X_i$ which is
isomorphic over the smooth locus of $\mathrm{Bl}_{Z_{i-1}}(X_{i-1})$.
We note that ${b'}_{\!i}$ is isomorphic at the generic point of
the center of $E$
on $\mathrm{Bl}_{Z_{i-1}}(X_{i-1})$.
We define $g_i = b_i \circ {b'}_{\!i} \colon X_i \to X_{i-1}$.
\item Define $Z_i$ as the center of $E$ on $X_i$ with the
reduced induced closed subscheme structure, and $F_i$ as
the only $g_i$-exceptional prime divisor on $X_i$ which contains $Z_i$.
\item We stop this process when $Z_n = F_n$.
This process must terminate after finite steps (see Remark \ref{finite})
and thus we get the sequence
$X_n \to \cdots \to X_0$.
\item We define $m \le n$ as the largest integer such that $Z_{m-1}$
is a point.
\item We define $g_{ji}\,(j \le i)$ as the morphism from $X_i$ to $X_j$.
\end{enumerate} 
\end{constr}

\begin{rem}\label{ideals}
We remark that $f_*\sO_Y(-iE) = g_{0n*}\sO_{X_n}(-iF_n)$ for any $i$
because $E$ and $F_n$ are the same as valuations.
\end{rem}

\begin{rem}\label{finite}
We prove the termination of the process.
Assume that we have the sequence $X_l \to \cdots \to X_0$ and $Z_l \neq F_l$.
We take common resolutions of $X_l$ and $Y$ over $X$, that is,
birational morphisms $h \colon W \to X_l$ and $h' \colon W \to Y$
from a smooth variety $W$ such that $g_{0l} \circ h = f \circ h'$.
We put
\begin{align*}
K_Y &= f^*K_X + aE, \\
K_{X_l} &= g_{0l}^*K_X + sF_l + (\mathrm{others}), \\
K_W &= h^*K_{X_l} + c({h'}^{-1})_*E + (\mathrm{others}), \\
h^*F_l &= (h^{-1})_*F_l + t({h'}^{-1})_*E + (\mathrm{others}).
\end{align*}
We note that $a, s, c$ and $t$ are positive integers. Then
\begin{align*}
K_W &= {h'}^*(f^*K_X + aE) + (\mathrm{others}) \\
&= h^*(g_{0l}^*K_X + sF_l + (\mathrm{others}))
+ c({h'}^{-1})_*E + (\mathrm{others}) \\
&= h^*g_{0l}^*K_X + s(h^{-1})_*F_l
+ (st + c)({h'}^{-1})_*E + (\mathrm{others}).
\end{align*}
Comparing the coefficients of $({h'}^{-1})_*E$, we have $a = st + c$
and especially $a > s$.
On the other hand because we know $s \ge l+1$ by the construction of $F_l$,
we get $a > l+1$.
It shows that the above process terminates with $n \le a-1$.
\qed
\end{rem}

We state an easy lemma.

\begin{lem}\label{excep}
Let  $f_i \colon (Y_i \supset E_i) \to (X \supset f_i (E_i))$ with $i = 1,2$
be algebraic germs of divisorial contractions.
Assume that $E_1$ and $E_2$ are the same as valuations.
Then $f_1$ and $f_2$ are isomorphic as morphisms over $X$.
\end{lem}

\begin{proof}
Let  $g_i \colon Z \to Y_i$ with $i = 1,2$ be common resolutions
and $h = f_i \circ g_i$.
We choose $g_i$-exceptional effective $\mQ$-divisors $F_i\,(i=1,2)$
and a $\mQ$-divisor $G$
on $Z$ such that $G = - g_1^*E_1 + F_1 = - g_2^*E_2 + F_2$. 
Then,
\begin{align*}
Y_i = \Proj_X\oplus_{j\ge0}f_{i*}\sO_{Y_i}(-jE_i)
= \Proj_X\oplus_{j\ge0}h_*\sO_Z(jG).
\end{align*}
\end{proof}

For weighted blow-ups in dimension $3$, we have a criterion
on terminal singularities.

\begin{thm}\label{toric}
Let $X \ni P$ be an algebraic germ
of a smooth $3$-dimensional variety with local parameters
$x, y, z$ at $P$, let $r, a, b$ be positive integers with $r \le a \le b$,
and let $Y \to X$ be the weighted blow-up of $X$ with its weights
$(x, y, z) = (r, a, b)$.
Then $Y$ has only terminal singularities if and only if $r = 1$ and
$a, b$ are coprime.
\end{thm}

By the above lemma and theorem, the problem is reduced to proving that
$F_n$ equals, as valuations, an exceptional divisor obtained by
a weighted blow-up of $X$.
We restate this in terms of ideal sheaves of $\sO_X$.

\begin{prop}\label{restate}
(Notation as above). 
$F_n$ equals, as valuations, an exceptional divisor obtained by
a weighted blow-up of $X$ with its weights $(x, y, z) = (1, m, n)$
for suitable local parameters $x, y, z$ at $P$,
if and only if the following conditions hold:

\smallskip
1. $f_*\sO_Y(-2E) \neq \gm_P$, that is, $g_{0n*}\sO_{X_n}(-2F_n) \neq \gm_P$.

2. $f_*\sO_Y(-nE) \not\subseteq \gm_P^2$,
that is, $g_{0n*}\sO_{X_n}(-nF_n) \not\subseteq \gm_P^2$.

\smallskip
\noindent Here $\gm_P \subset \sO_X$ is the ideal sheaf of $P$.
\end{prop}

\begin{proof}
The ``only if'' part is obvious taking it into account that for any $i$
$g_{0n*}\sO_{X_n}(-iF_n) = (x^sy^tz^u|s+mt+nu \ge i)$.
Actually $x \not\in g_{0n*}\sO_{X_n}(-2F_n)$ and
$z \in g_{0n*}\sO_{X_n}(-nF_n)$.

Now we prove the ``if'' part.
The condition 1 means that the coefficient of $F_n$ in $g_{1n}^*F_1$ is $1$.
This says that for any $i \ge 1$, $F_i$ is the only $g_{0i}$-exceptional
prime divisor on $X_i$ containing $Z_i$
and the coefficient of $F_n$ in $g_{in}^*F_i$ is $1$.

We consider a prime divisor $D \ni P$ on $X$ which is smooth at $P$
and define $1 \le l \le n$ as the largest integer such that
$Z_{l-1} \subseteq (g_{0,l-1}^{-1})_*D$.
Then $(g_{0i}^{-1})_*D$ is smooth at the generic
point of $Z_i$ for any $i<l$, and so we get
$g_{0l}^*D = (g_{0l}^{-1})_*D +\sum_{i=1}^li\,(g_{il}^{-1})_*F_i
+(\mathrm{others})$. 
Therefore the coefficient of $F_n$ in $g_{0n}^*D$ is $l$.
By the condition 2, we can choose $z \in \gm_P\!\setminus\!\gm_P^2$
such that $g_{0n}^*\wdiv(z) \ge nF_n$, that is,
$Z_{n-1} \subseteq (g_{0,n-1}^{-1})_*\wdiv(z)$
because of the above argument.
Adding $x, y \in \gm_P\!\setminus\!\gm_P^2$ such that
$Z_{m-1} \subseteq (g_{0,m-1}^{-1})_*\wdiv(y)$,
we can take local parameters $x, y, z$ at $P$. 
Then $F_i\,(1 \le i \le n)$ equals, as valuations,
the exceptional divisor obtained by the weighted blow-up of $X$
with its weights $(x, y, z) = (1, \min\{i, m\}, i)$,
and especially $F_n$ is obtained by the weighted blow-up of $X$
with its weights $(x, y, z) = (1, m, n)$.
\end{proof}

So we prove the above two conditions.

\subsection{Preliminaries}
Let $K_Y = f^*K_X + aE$, and let $r$ be the global
Gorenstein index of $Y$, that is, the smallest positive integer
such that $rK_Y$ is Cartier. 
Since $a$ equals the discrepancy of $F_n$
with respect to $K_X$, $a \in \mZ_{\ge 2}$.

\begin{lem}\label{coprime}
(Notation as above). $a$ and $r$ are coprime.
\end{lem}

\begin{proof}
Let $s$ be the greatest common divisor of $a$ and $r$,
and let $a = sa', r = sr'.$
Since $r'aE = a'rE$ is Cartier by \cite[Corollary 5.2]{K88}, so is $r'K_Y$.
Hence $r'=r$ and $s=1$.
\end{proof}

We recall the singular Riemann-Roch formula (\cite[Theorem 10.2]{R87}).

\begin{thm}\label{singRR}
Let $X$  be a projective $3$-dimensional variety
with only canonical singularities,
and let $D$ be a Weil divisor on $X$ such that for any $P \in X$
there exists an integer $i_P$
satisfying $(\sO_X(D))_P \cong (\sO_X(i_PK_X))_P$.
Then there is a formula of the form
\begin{align*}
\chi(\sO_X(D)) &= \chi(\sO_X) + \frac{1}{12} D(D - K_X)(2D- K_X) \\
&\quad + \frac{1}{12} D\cdot c_2(X) + \sum_Pc_P(D),
\end{align*}
where the summation takes place over singular points of $X$,
and $c_P(D) \in \mQ$ is a contribution depending only on the local analytic
type of $P \in X$ and $\sO_X(D)$.

If $P$ is a terminal quotient singularity of
type  $\frac{1}{r_P}(1, -1, b_P)$, then
\begin{align*}
c_P(D) = - \overline{i_P} \frac{r_P^2-1}{12r_P} +
\sum_{j=1}^{\overline{i_P}-1}
\frac{\overline{jb_P}(r_P-\overline{jb_P})}{2r_P},
\end{align*}
where $\bar{\ }$ denotes the smallest residue modulo $r_P$,
that is, $\overline{j} = j - \lfloor \frac{j}{r_P} \rfloor r_P$
in terms of the round down $\lfloor \ \rfloor$.
The definition of the round down $\lfloor \ \rfloor$ is
$\lfloor j \rfloor = \max\{k \in \mZ | k \le j\}$.

And for any terminal singularity $P$, 
\begin{align*}
c_P(D) = \sum_{\alpha} c_{P_{\alpha}}(D_{\alpha}),
\end{align*}
where $\{(P_{\alpha},D_{\alpha})\}_{\alpha}$ is a flat deformation of $(P, D)$
to terminal quotient singularities.
\end{thm}

\begin{rem}\label{singRRrem}
If $X$ has only terminal singularities, then we can write
the contribution term $\sum_Pc_P(D)$ as $\sum_Qc_Q(D)$,
where
\begin{align*}
c_Q(D) = - \overline{i_Q} \frac{r_Q^2-1}{12r_Q} +
\sum_{j=1}^{\overline{i_Q}-1}
\frac{\overline{jb_Q}(r_Q-\overline{jb_Q})}{2r_Q}.
\end{align*}
For its summation takes place over points which need not lie on $X$
but may lie on deformed varieties of $X$, $Q$'s are called ``fictitious''
points in the sense of M.\@ Reid.
This description holds even though $X$ has canonical singularities,
but in this case $Q$'s may lie on deformed varieties of crepant
blown-up varieties of $X$ (see \cite{R87} for details).
\end{rem}

By Lemma \ref{coprime}, we can take an integer $e$ such that
$ae \equiv 1$ modulo $r$. 
Then $(\sO_Y(E))_Q \cong (\sO_Y(eK_Y))_Q$ for any $Q \in E$.
Using the singular Riemann-Roch formula, we get
\begin{align}
\chi(\sO_Y(iE)) &= \chi(\sO_Y) + \frac{1}{12} i(i - a)(2i - a)E^3
\label{eqRR} \\
\nonumber &\quad + \frac{1}{12} i E\cdot c_2(Y) + A_i,
\end{align}
where $A_i$ is the contribution term and has the below description:
\begin{gather*}
A_i =\sum_{Q\in I}c_Q(iE), \\
c_Q(iE) = - \overline{ie} \frac{r_Q^2-1}{12r_Q} +
\sum_{j=1}^{\overline{ie}-1}
\frac{\overline{jb_Q}(r_Q-\overline{jb_Q})}{2r_Q}.
\end{gather*}
Here $Q\in I$ are fictitious singularities. 
The type of $Q$ is $\frac{1}{r_Q}(1, -1, b_Q)$,
$(\sO_{Y_Q}(E_Q))_Q \cong (\sO_{Y_Q}(eK_{Y_Q}))_Q$
where $(Y_Q, E_Q)$ is the fictitious pair for $Q$,
and $\bar{\ }$ denotes the smallest residue modulo $r_Q$.
We note that $b_Q$ is coprime to $r_Q$ and also
$e$ is coprime to $r_Q$ because $r|(ae-1)$.
So $v_Q = \overline{eb_Q}$ is coprime to $r_Q$.
With this description, $r = 1$ if $I$ is empty, and otherwise
$r$ is the lowest common multiple of $\{r_Q\}_{Q \in I}$.
We note that $c_Q(iE)$ depends only on $i$ mod $r_Q$ and
equals $0$ if $r_Q|i$. Especially $A_i$ depends only on $i$ mod $r$
and equals $0$ if $r|i$.

We put $B_i = -(A_i + A_{-i})$.
Because
\begin{align*}
c_Q(iE) + c_Q(-iE)
&= \Bigl( -\overline{ie} \frac{r_Q^2-1}{12r_Q} 
+ \sum_{j=1}^{\overline{ie}-1}
\frac{\overline{jb_Q}(r_Q-\overline{jb_Q})}{2r_Q} \Bigr) \\
&\qquad + \Bigl( -\overline{-ie} \frac{r_Q^2-1}{12r_Q} 
+ \sum_{j=1}^{\overline{-ie}-1}
\frac{\overline{jb_Q}(r_Q-\overline{jb_Q})}{2r_Q} \Bigr) \\
&= -\frac{r_Q^2-1}{12} + 
\Bigl( \sum_{j=1}^{r_Q}
\frac{\overline{jb_Q}(r_Q-\overline{jb_Q})}{2r_Q} \Bigr)
-  \frac{\overline{ieb_Q}(r_Q-\overline{ieb_Q})}{2r_Q} \\
&= -\frac{r_Q^2-1}{12} + 
\Bigl( \sum_{j=1}^{r_Q}
\frac{j(r_Q-j)}{2r_Q} \Bigr)
-  \frac{\overline{iv_Q}(r_Q-\overline{iv_Q})}{2r_Q} \\
&= - \frac{\overline{iv_Q}(r_Q-\overline{iv_Q})}{2r_Q}
\end{align*}
where the third equality comes from
the property that $b_Q$ and $r_Q$ are coprime,
we have
\begin{align}
B_i = - \sum_{Q\in I}(c_Q(iE) + c_Q(-iE)) 
= \sum_{Q\in I} \frac{\overline{iv_Q}(r_Q-\overline{iv_Q})}{2r_Q}. \label{Bi}
\end{align}

\begin{prop}\label{eqs}
(Notation as above).
\begin{align}
&rE^3 \in \mZ_{>0}. \tag{A}\label{AA} \\
&1 = \frac{1}{2} aE^3 + \sum_{Q\in I} \frac{v_Q(r_Q - v_Q)}{2r_Q}.
\tag{B}\label{BB} \\
&\dim_k\sO_X/f_*\sO_Y(-iE) = \tag{C}\label{CC} \\
&\qquad i^2 - \frac{1}{2}\sum_{Q \in I}
\min_{0 \le j < i}\{(1+j)jr_Q + i(i-1-2j)v_Q\} \quad (1 \le i \le a).
\notag \\
&\sum_{Q \in I} \min \{v_Q, r_Q - v_Q\} = \dim_kf_*\sO_Y(-2E)/\gm_P^2.
\tag{D}\label{DD}
\end{align}
\end{prop}

\begin{rem}\label{remeqs}
In particular (\ref{AA}), (\ref{CC}) and (\ref{DD}) are essential.
We use (\ref{AA}) to bound the value of $a$ from above and
use (\ref{CC}) to control the values of $r_Q$'s.
(\ref{DD}) shows that the number of fictitious non-Gorenstein points of $Y$
is at most $3$. We prove the conditions 1 and 2 in Proposition \ref{restate}
according to the value of $\dim_kf_*\sO_Y(-2E)/\gm_P^2$.
\end{rem}

\begin{rem}\label{remCC}
In fact, because of (\ref{Bi}) and (\ref{sumchi2-2})
the right hand side of (\ref{CC}) is the same if we replace
$v_Q$ by $r_Q - v_Q$.
\end{rem}

\begin{proof}
We consider the exact sequence:
\begin{align}
0 \to \sO_Y((i-1)E) \to \sO_Y(iE) \to \sQ_i \to 0. \label{exact}
\end{align}
By (\ref{eqRR}), we get
\begin{align}
\chi(\sQ_i) &= \chi(\sO_Y(iE)) - \chi(\sO_Y((i-1)E)) \label{Q} \\
\nonumber &= \frac{1}{12} \{2(3i^2-3i+1)
- 3(2i-1)a + a^2\}E^3 \\
\nonumber &\quad + \frac{1}{12} E\cdot c_2(Y) + A_i - A_{i-1}.
\end{align}
Since $\chi(\sQ_i) - \chi(\sQ_{r+i})
= \frac{r}{2}(a+1-r-2i)E^3$ is an integer for
any $i$ and $E^3$ is positive, we have (\ref{AA}).

By (\ref{Q}),
\begin{align}
\chi(\sQ_{-i}) - \chi(\sQ_{i+1}) = (i+\frac{1}{2})aE^3 + B_{i+1} - B_i.
\label{diffchi}
\end{align}

Let $d(i) = \dim_kf_*\sO_Y(iE)/f_*\sO_Y((i-1)E)$.
We note that $d(i) = 0$ if $i \ge 1$, and $d(0) = 1$.
Because $(Y, \varepsilon E)$ is weak KLT and $iE - (K_Y + \varepsilon E)$
is $f$-ample for a sufficiently small positive rational number $\varepsilon$
and an integer $i \le a$, using \cite[Theorem 1-2-5]{KMM87},
we have $R^jf_*\sO_Y(iE) = 0$ for $i \le a,\, j \ge 1$.
So by (\ref{exact}), for any $i \le a$,
\begin{itemize}
\item $H^0(Y, \sQ_i) = f_*\sQ_i = f_*\sO_Y(iE)/f_*\sO_Y((i-1)E),$  
\item $H^j(Y, \sQ_i) = R^jf_*\sQ_i = 0$ \quad for \ $j \ge 1$,
\end{itemize}
and therefore $d(i) = \chi(\sQ_i)$.

Putting $i=0$ in (\ref{diffchi}), we get
\begin{align}
&1 = \frac{1}{2}aE^3 + B_1. \label{B1}
\end{align}
Combining this and (\ref{Bi}) with $i=1$, we get (\ref{BB}).

With (\ref{diffchi}), we obtain for $1 \le i \le a$,
\begin{align}
\sum_{1 \le j < i}d(-j) &=
\sum_{1 \le j < i}\{\chi(\sQ_{-j}) - \chi(\sQ_{j+1})\} \label{sumchi} \\
\nonumber &= \sum_{1 \le j < i}\{(j+\frac{1}{2})aE^3 + B_{j+1} - B_j\} \\
\nonumber &= \frac{1}{2}(i^2-1)aE^3 + B_i - B_1.
\end{align}
Eliminating $\frac{1}{2} aE^3$ with (\ref{B1}), we obtain
\begin{align}
\sum_{1 \le j < i}d(-j) = (i^2-1) + B_i - i^2 B_1
\quad(1 \le i \le a). \label{sumchi2-1}
\end{align}
Since for $i \ge 1$,
\begin{align*}
&\quad \frac{\overline{iv_Q}(r_Q-\overline{iv_Q})}{2r_Q}
- i^2 \frac{v_Q(r_Q-v_Q)}{2r_Q} \\
&= -\frac{1}{2} 
\Bigm\{ r_Q\Bigl( \frac{iv_Q - \overline{iv_Q}}{r_Q}
- \frac{iv_Q}{r_Q} + \frac{1}{2} \Bigr)^2 + i^2
\frac{v_Q(r_Q-v_Q)}{r_Q} - \frac{r_Q}{4} \Bigm\} \\
&= -\frac{1}{2}
\Bigm\{ r_Q\Bigl(\Bigm\lfloor \! \frac{iv_Q}{r_Q} \! \Bigm\rfloor
- \frac{iv_Q}{r_Q} + \frac{1}{2} \Bigr)^2 + i^2\frac{v_Q(r_Q-v_Q)}{r_Q}
- \frac{r_Q}{4} \Bigm\} \\
&= -\frac{1}{2} \min_{0 \le j < i}
\Bigm\{ r_Q\Bigl( j - \frac{iv_Q}{r_Q} + \frac{1}{2} \Bigr)^2
+ i^2\frac{v_Q(r_Q-v_Q)}{r_Q} - \frac{r_Q}{4} \Bigm\} \\
&= -\frac{1}{2} \min_{0 \le j < i} \{(1+j)jr_Q + i(i-1-2j)v_Q\},
\end{align*}
with (\ref{Bi}) we have
\begin{align}
B_i - i^2 B_1 \label{sumchi2-2}
= - \frac{1}{2}\sum_{Q \in I}\min_{0 \le j < i}\{(1+j)jr_Q + i(i-1-2j)v_Q\}
\quad (i \ge 1).
\end{align}
Of course because $\dim_k\sO_X/f_*\sO_Y(-iE) = 1 + \sum_{j=1}^{i-1}d(-j)$,
combining this with (\ref{sumchi2-1}) and (\ref{sumchi2-2}),
we obtain (\ref{CC}).

Putting $i = 2$ in (\ref{sumchi2-1}) and (\ref{sumchi2-2}), we have
\begin{align*}
d(-1) = 3 - \sum_{Q \in I} \min \{v_Q, r_Q - v_Q\}.
\end{align*}
Since $\dim_kf_*\sO_Y(-2E)/\gm_P^2 = 3 - d(-1)$, we get (\ref{DD}).
\end{proof}

\subsection{Proof of $\mathbf{f_*\sO_Y(-2E) \neq \gm_P}$}\label{2E}
Assuming that $f_*\sO_Y(-2E) = \gm_P$, we will derive a contradiction.
The assumption means that the coefficient of $F_n$ in
$g_{1n}^*F_1$ is bigger than $1$, so there exists a $Z_i$
which is contained in at least two $g_{0i}$-exceptional
prime divisors on $X_i$. 
The minimum value of $a$ in this case occurs when
$Z_1$ is a curve, $Z_2 = (g_2^{-1})_*F_1 \cap F_2$, and $n=3$,
and the minimum value is $6$. So we get $a \ge 6$.
By the assumption and (\ref{DD}),
we obtain $\sum_{Q \in I} \min \{v_Q, r_Q - v_Q\} = 3$.
Thus we have only to consider the three cases:

\smallskip
\begin{tabular}{l}
Case\,1. $\{(r_Q, \overline{v_Q})\}_{Q \in I} =
\{(r, \overline{\pm3})\}, \ r \ge 7.$ \\
Case\,2. $\{(r_Q, \overline{v_Q})\}_{Q \in I} =
\{(r_1, \overline{\pm1}), (r_2, \overline{\pm2})\},
\ r_1 \ge 2, \ r_2 \ge 5.$ \\
Case\,3. $\{(r_Q, \overline{v_Q})\}_{Q \in I} =
\{(r_1, \overline{\pm1}), (r_2, \overline{\pm1}), (r_3, \overline{\pm1})\},
\ 2 \le r_1 \le r_2 \le r_3.$ \\
\end{tabular}

\smallskip
\noindent Here $\pm$ means that one of these occurs for each $\overline{v_Q}$.
We remark that $v_Q$ is coprime to $r_Q$.

Since
$\displaystyle \sum_{Q\in I} \frac{v_Q(r_Q - v_Q)}{2r_Q} < 1$
from (\ref{BB}), we have the below inequalities:

\smallskip
\begin{tabular}{l}
Case\,1. $3/2 - 9/2r <1$. \\
Case\,2. $3/2 - (1/2r_1 + 2/r_2) <1$. \\
Case\,3. $3/2 - (1/2r_1 + 1/2r_2 + 1/2r_3) <1$.
\end{tabular}

\smallskip
Using this evaluation,
we can restrict possible values of $r_Q$'s.
Below we show all the possible values and the corresponding values of $aE^3$:

\smallskip
\begin{tabular}{ccccc}
Case\,1. & $r$ &:& $7$ & $8$ \\
 & $aE^3$ &:& $2/7$ & $1/8$
\end{tabular}

\begin{tabular}{ccccccc}
Case\,2. & $(r_1, r_2)$ &:& $(2, 5)$ & $(3, 5)$ & $(4, 5)$ & $(2, 7)$ \\
 & $aE^3$ &:& $3/10$ & $2/15$ & $1/20$ & $1/14$
\end{tabular}

\begin{tabular}{ccccccc}
Case\,3. & $(r_1, r_2, r_3)$ &:& $(2, 2, r_3)$ & $(2, 3, 3)$
& $(2, 3, 4)$ & $(2, 3, 5)$ \\
 & $aE^3$ &:& $2/2r_3$ & $1/6$ & $1/12$ & $1/30$
\end{tabular}

\smallskip
Recalling that $r$ is the lowest common multiple of $\{r_Q\}_{Q \in I}$,
with (\ref{AA}) we have $a \le 3$ for all the above cases.
This contradicts $a \ge 6$. 
\qed

\subsection{Proof of $\mathbf{f_*\sO_Y(-nE) \not\subseteq \gm_P^2}$}\label{nE}
Because $f_*\sO_Y(-2E) \neq \gm_P$, we have
$g_{1n}^*F_1 = \sum_{i=1}^n (g_{in}^{-1})_*F_i+(\mathrm{others})$ and,
\renewcommand{\labelenumi}{$(*)$}
\begin{enumerate}
\item
$F_i\,(1 \le i \le m)$ is obtained as a valuation by the weighted blow-up
of $X$ with its weights $(x, y, z) = (1, i, i)$ for local parameters
$x, y, z$ at $P$ such that $Z_{m-1} \subseteq
(g_{0,m-1}^{-1})_*\wdiv(y) \cap (g_{0, m-1}^{-1})_*\wdiv(z)$.
\end{enumerate}
We divide the proof according to the value of
$\dim_kf_*\sO_Y(-2E)/\gm_P^2 \le 2$.

\smallskip
\begin{tabular}{ll}
Case\,1.& $\dim_kf_*\sO_Y(-2E)/\gm_P^2 = 0$. \\
& This is the case when $Z_1 \subseteq F_1$ is
neither a line nor a point. \\
Case\,2.& $\dim_kf_*\sO_Y(-2E)/\gm_P^2 = 1$. \\
& This is the case when $Z_1 \subseteq F_1$ is a line. \\
Case\,3.& $\dim_kf_*\sO_Y(-2E)/\gm_P^2 = 2$. \\
& This is the case when $Z_1 \subseteq F_1$ is a point.
\end{tabular}

\smallskip
\noindent Since
\begin{align*}
\dim_kf_*\sO_Y(-2E)/\gm_P^2
&= \dim_k \mathrm{Im}[(v \in \gm_P | Z_1 \subseteq (g_1^{-1})_*\wdiv(v))
\to \gm_P/\gm_P^2] \\
&= \dim_k\{v \in \varGamma(F_1, \sO_{F_1}(1)) |
v = 0 \ \mathrm{or} \ Z_1 \subseteq \wdiv(v)\},
\end{align*}
the value of $\dim_kf_*\sO_Y(-2E)/\gm_P^2$
decides the type of $Z_1 \subseteq F_1 \cong \mP_k^2$
as above.

In Case\,1, $\sum_{Q \in I} \min \{v_Q, r_Q - v_Q\} = 0$ by (\ref{DD}).
Therefore $I$ is empty and thus $Y$ is Gorenstein.
By \cite[Theorem 5]{Cu88}, $f$ must be the blow-up of $X$ along $P$,
that is, $f = g_1$, and so we have nothing to do.
Thus we have only to consider Cases\,\,2 and 3.
In these cases we investigate the values of $\dim_k\sO_X/f_*\sO_Y(-iE)$'s
carefully.

\begin{prop}\label{keyprop}
(Notation as above).
Let $2 \le l \le n$ be an integer such that
$g_{0i*}(-iF_i) \not\subseteq \gm_P^2$ for any $i<l$.

(1) If $g_{0l*}(-lF_l) \not\subseteq \gm_P^2$, then
\begin{align*}
\dim_k\sO_X/f_*\sO_Y(-lE)
\le l - \frac{1}{2} \min_{0 \le j < l}\{((1+j)m-2l)j\}.
\end{align*}

(2) If $g_{0l*}(-lF_l) \subseteq \gm_P^2$
(in this case we have $l > m$ by $(*)$), then
\begin{align*}
\dim_k\sO_X/f_*\sO_Y(-lE)
> l - \frac{1}{2} \min_{0 \le j < l}\{((1+j)m-2l)j\}.
\end{align*}
\end{prop}

\begin{rem}\label{keyrem}
In the case when $m = 1$ because
\begin{align*}
\min_{0 \le j < l}\{((1+j)m-2l)j\} = \min_{0 \le j < l}\{(j-(2l-1))j\}
= -l(l-1),
\end{align*}
we can simplify the above inequalities:

\smallskip
(1) $\displaystyle \dim_k\sO_X/f_*\sO_Y(-lE) \le \frac{1}{2}l(l+1).$

\smallskip
(2) $\displaystyle \dim_k\sO_X/f_*\sO_Y(-lE) > \frac{1}{2}l(l+1).$
\end{rem}

\begin{proof}
\noindent (1) By the assumption and $f_*\sO_Y(-2E) \neq \gm_P$,
the proof of Proposition \ref{restate} says that we can take local parameters
$x, y, z$ at $P$ such that
$Z_{\min\{l, m\}-1} \subseteq (g_{0,\min\{l, m\}-1}^{-1})_*\wdiv(y)$
and $Z_{l-1} \subseteq (g_{0,l-1}^{-1})_*\wdiv(z)$.
Then for $1 \le i \le l$, $F_i$ equals, as valuations,
the exceptional divisor obtained by the weighted blow-up of $X$
with its weights $(x, y, z) = (1, \min\{i, m\}, i)$.

Hence
\begin{align*}
&\quad f_*\sO_Y(-lE) = g_{0n*}\sO_{X_n}(-lF_n) \\
\nonumber &\supseteq g_{0l*}\sO_{X_l}(-lF_l)
= (x^sy^tz^u|s+\min\{l, m\}t+lu \ge l),
\end{align*}
and so
\begin{align*}
\dim_k\sO_X/f_*\sO_Y(-lE) &\le
\dim_k\sO_X/(x^sy^tz^u|s+\min\{l, m\}t+lu \ge l) \\
&= l - \frac{1}{2} \min_{0 \le j < l}\{((1+j)m-2l)j\}.
\end{align*}
Here we used Lemma \ref{spandim} proved later.

(2) As in the proof of (1), we can take local parameters
$x, y, z$ at $P$ such that $Z_{m-1} \subseteq (g_{0,m-1}^{-1})_*\wdiv(y)$
and $Z_{l-2} \subseteq (g_{0,l-2}^{-1})_*\wdiv(z)$.
Then for $1 \le i < l$, $F_i$ equals, as valuations,
the exceptional divisor obtained by the weighted blow-up of $X$
with its weights $(x, y, z) = (1, \min\{i, m\}, i)$.

We have
\begin{align*}
&\quad f_*\sO_Y(-lE) = g_{0n*}\sO_{X_n}(-lF_n) \\
&\subseteq g_{0,l-1*}\sO_{X_{l-1}}(-lF_{l-1})
+ (v \in \gm_P | Z_{l-1} \subseteq (g_{0,l-1}^{-1})_*\wdiv(v)).
\end{align*}
But since
\begin{align*}
(v \in \gm_P | Z_{l-1} \subseteq (g_{0,l-1}^{-1})_*\wdiv(v))
\subseteq g_{0l*}\sO_{X_l}(-lF_l) \subseteq \gm_P^2,
\end{align*}
for any $v \in \gm_P$
such that $Z_{l-1} \subseteq (g_{0,l-1}^{-1})_*\wdiv(v)$ we have
\begin{align*}
g_{0,l-1}^*\wdiv(v) \ge g_{1,l-1}^*(2F_1 + (g_1^{-1})_*\wdiv(v))
\ge (2+(l-2))F_{l-1} = lF_{l-1}.
\end{align*}
Thus
\begin{align*}
(v \in \gm_P | Z_{l-1} \subseteq (g_{0,l-1}^{-1})_*\wdiv(v))
\subseteq g_{0,l-1*}\sO_{X_{l-1}}(-lF_{l-1}),
\end{align*}
and hence
\begin{align*}
f_*\sO_Y(-lE)  \subseteq g_{0,l-1*}\sO_{X_{l-1}}(-lF_{l-1})
= (x^sy^tz^u|s+mt+(l-1)u \ge l).
\end{align*}
Therefore with Lemma \ref{spandim},
\begin{align*}
\dim_k\sO_X/f_*\sO_Y(-lE) &\ge \dim_k\sO_X/(x^sy^tz^u|s+mt+(l-1)u \ge l) \\
&> \dim_k\sO_X/(x^sy^tz^u|s+mt+lu \ge l) \\
&= l - \frac{1}{2} \min_{0 \le j < l}\{((1+j)m-2l)j\}.
\end{align*}
\end{proof}

We used the following lemma in the above proof.

\begin{lem} \label{spandim}
Let $X \ni P$ be an algebraic germ of a smooth $3$-dimensional variety
with local parameters $x, y, z$ at $P$, and let $l, m$ be positive integers.
Then
\begin{align*}
\dim_k\sO_X/(x^sy^tz^u|s+\min\{l, m\}t+lu \ge l)
= l - \frac{1}{2} \min_{0 \le j < l}\{((1+j)m-2l)j\}.
\end{align*}
\end{lem}

\begin{proof}
\begin{align*}
&\quad \dim_k\sO_X/(x^sy^tz^u|s+\min\{l, m\}t+lu \ge l) \\
&= \dim_k\Span_k\langle x^sy^t|s+\min\{l, m\}t<l\rangle \\
&= \sum_{0 \le t \le \lfloor \frac{l}{\min\{l, m\}} \rfloor}
(l - \min\{l, m\}t) \\
&= \sum_{0 \le t \le \lfloor \frac{l}{m} \rfloor} (l - mt) \\
&= l - \frac{m}{2} \Bigm\{
\Bigl( \Bigm\lfloor \! \frac{l}{m} \! \Bigm\rfloor
+ \frac{1}{2} - \frac{l}{m} \Bigr)^2
- \Bigl( \frac{1}{2} - \frac{l}{m} \Bigr)^2 \Bigm\} \\
&= l - \frac{m}{2}\min_{0 \le j < l}
\Bigm\{ \Bigl( j + \frac{1}{2} - \frac{l}{m} \Bigr)^2
- \Bigl( \frac{1}{2} - \frac{l}{m} \Bigr)^2 \Bigm\} \\
&= l - \frac{1}{2} \min_{0 \le j < l}\{((1+j)m-2l)j\}.
\end{align*}
\end{proof}

Now we prove $f_*\sO_Y(-nE) \not\subseteq \gm_P^2$ in Cases\,\,2 and 3.

\begin{proof}[Proof in Case\,2]
For $Z_1 \subseteq F_1$ is a line in this case
and $a$ is the discrepancy of $F_n$ with respect to $K_X$, we get $m=1$ and
\begin{align}
a = n+1 \quad (n \ge 2). \label{a}
\end{align}
By (\ref{DD}), $\sum_{Q \in I} \min \{v_Q, r_Q - v_Q\} = 1$
and thus $\{(r_Q, \overline{v_Q})\}_{Q \in I} = \{(r, \overline{\pm1})\}$.
From (\ref{BB}), we obtain $aE^3 = (r + 1)/r$.
By (\ref{AA}),
\begin{align}
a \le r+1. \label{ar}
\end{align}
From (\ref{CC}), Remark \ref{remCC}, and (\ref{a}),
for $1 \le i \le n+1$ we have
\begin{align*}
\dim_k\sO_X/f_*\sO_Y(-iE)
&= i^2 - \frac{1}{2}\min_{0 \le j < i}\{(1+j)jr + i(i-1-2j)\} \\
&= \frac{1}{2}i(i+1)-\frac{1}{2}\min_{0 \le j < i}\{((1+j)r-2i)j\}.
\end{align*}
Hence for $1 \le i \le n+1$,
\begin{align}
\dim_k\sO_X/f_*\sO_Y(-iE) \ge \frac{1}{2}i(i+1), \label{dim}
\end{align}
where the equality holds if and only if $i \le r$.

If there exists a positive integer $2 \le l \le n$
such that $g_{0l*}\sO_{X_l}(-lF_l) \subseteq \gm_P^2$ and
$g_{0i*}\sO_{X_i}(-iF_i) \not\subseteq \gm_P^2$ for any $i<l$,
then by Proposition \ref{keyprop}, Remark \ref{keyrem}, and
the condition of the equality in (\ref{dim}), we obtain $l = r+1$.
Thus with (\ref{a}),  we have $r+1 = l \le n = a-1$, that is, $a \ge r+2$.
This contradicts (\ref{ar}) and hence we get
$g_{0n*}\sO_{X_n}(-nF_n) \not\subseteq \gm_P^2$.
\end{proof}

\begin{proof}[Proof in Case\,3]
In this case we use essentially
the same idea as in Case\,2, but it is a little more complicated.
By (\ref{DD}),
$\sum_{Q \in I} \min \{v_Q, r_Q - v_Q\} = 2$.
Thus we have only to consider the two subcases:

\begin{itemize}
\item Subcase\,1. $\{(r_Q, \overline{v_Q})\}_{Q \in I} =
\{(r, \overline{\pm2})\}, \ r \ge 5.$
\item Subcase\,2. $\{(r_Q, \overline{v_Q})\}_{Q \in I} =
\{(r_1, \overline{\pm1}), (r_2, \overline{\pm1})\},
\ 2 \le r_1 \le r_2.$
\end{itemize}

In Subcase\,1, we have $aE^3 = 4/r$ by (\ref{BB})
and thus $a \le 4$ from (\ref{AA}).
But since $Z_1 \subseteq F_1$ is a point, we get
$n = 2$ and $a = 4$.
Then choosing local parameters $x, y, z$ at $P$
such that $Z_1 \subseteq (g_1^{-1})_*\wdiv(y) \cap (g_1^{-1})_*\wdiv(z)$,
$F_2$ equals, as valuations, the exceptional divisor
obtained by the weighted blow-up of $X$
with its weights $(x, y, z) = (1, 2, 2)$.
So we have only to investigate Subcase\,2.

Recalling that $a$ is the discrepancy of $F_n$ with respect to $K_X$,
we have
\begin{align}
a = m + n \quad (2 \le m \le n). \label{a2}
\end{align}
Calculating with (\ref{BB}) we obtain $aE^3 = (r_1 + r_2)/r_1r_2$,
and thus by (\ref{AA}),
\begin{align}
a \le r_1+r_2. \label{ar2}
\end{align}
From (\ref{CC}), Remark \ref{remCC}, and (\ref{a2}),
for $1 \le i \le m+n$ we have
\begin{align*}
&\quad \dim_k\sO_X/f_*\sO_Y(-iE) \\
&= i^2 - \frac{1}{2}\min_{0 \le j < i} \{(1+j)jr_1 + i(i-1-2j)\} \\
&\qquad - \frac{1}{2}\min_{0 \le j < i} \{(1+j)jr_2 + i(i-1-2j)\} \\
&= i - \frac{1}{2} \Bigl( \min_{0 \le j < i}\{((1+j)r_1-2i)j\}
+ \min_{0 \le j < i}\{((1+j)r_2-2i)j\} \Bigr).
\end{align*}
Hence for $1 \le i \le m+n$,
\begin{align}
\dim_k\sO_X/f_*\sO_Y(-iE)
&\ge i - \frac{1}{2} \min_{0 \le j < i}\{((1+j)r_1-2i)j\} \label{dim2} \\
\nonumber &\ge i,
\end{align}
where the equality of the first inequality holds
if and only if $i \le r_2$, and the second holds if and only if $i \le r_1$.

\begin{cl}\label{r1}
$r_1 = m$.
\end{cl}

\begin{proof}[Proof of the claim]
Utilizing Proposition \ref{keyprop}\,(1) with $l=m$, we have
\begin{align}
\dim_k\sO_X/f_*\sO_Y(-mE)
\le m - \frac{1}{2} \min_{0 \le j < m}\{(j-1)jm\} = m. \label{dim3}
\end{align}
We take local parameters $x, y, z$ at $P$ as in $(*)$,
satisfying $Z_m \subseteq (g_{0m}^{-1})_*\wdiv(z)$
if $Z_m \subseteq F_m \cong \mP_k^2$ is a line.
We have
\begin{align*}
&\quad f_*\sO_Y(-(m+1)E)=g_{0n*}\sO_{X_n}(-(m+1)F_n) \\
\nonumber &\subseteq g_{0m*}\sO_{X_m}(-(m+1)F_m)
+(v\in \gm_P|Z_m \subseteq (g_{0m}^{-1})_*\wdiv(v)).
\end{align*}
But since
\begin{align*}
(v\in \gm_P|Z_m \subseteq (g_{0m}^{-1})_*\wdiv(v))
\subseteq g_{0m*}\sO_{X_m}(-mF_m)=(x^m, y, z),
\end{align*}
we get
\begin{align*}
(v\in \gm_P|Z_m \subseteq (g_{0m}^{-1})_*\wdiv(v))
\subseteq (z) + g_{0m*}\sO_{X_m}(-(m+1)F_m),
\end{align*}
and thus
\begin{align*}
f_*\sO_Y(-(m+1)E) &\subseteq (z) + g_{0m*}\sO_{X_m}(-(m+1)F_m) \\
&= (z) + (x^sy^tz^u|s+mt+mu \ge m+1).
\end{align*}
Hence,
\begin{align}
&\quad \dim_k\sO_X/f_*\sO_Y(-(m+1)E) \label{dim4} \\
\nonumber &\ge \dim_k\sO_X/((z) + (x^sy^tz^u|s+mt+mu \ge m+1)) \\
\nonumber &= \dim_k\Span_k\langle x^s,\,y | s \le m \rangle \\
\nonumber &= m + 2.
\end{align}
From (\ref{dim3}), (\ref{dim4}), and
the condition of the second equality in (\ref{dim2}), we have $r_1 = m$.
\end{proof}

If there exists a positive integer $l \le n$
such that $g_{0l*}\sO_{X_l}(-lF_l) \subseteq \gm_P^2$ and
$g_{0i*}\sO_{X_i}(-iF_i) \not\subseteq \gm_P^2$ for any $i<l$,
then by Proposition \ref{keyprop}, Claim \ref{r1}, and
the condition of the first equality in (\ref{dim2}), we obtain $l = r_2+1$.
Thus with (\ref{a2}) and Claim \ref{r1}, we have $r_1+r_2+1 = m+l \le m+n =a$.
This contradicts (\ref{ar2}) and hence we get
$g_{0n*}\sO_{X_n}(-nF_n) \not\subseteq \gm_P^2$.
\end{proof}

\bibliographystyle{amspalpha}

\begin{thebibliography}{KMM87}
\bibitem[Co95]{Co95}A.\@ Corti,
{\itshape Factoring birational maps of threefolds after Sarkisov},
Jour. Alg. Geom. 4 (1995), 223-254.
\bibitem[Co99]{Co99}A.\@ Corti,
{\itshape Singularities of linear systems and $3$-fold birational geometry},
Warwick, preprint, (1999).
\bibitem[CPR99]{CPR99}A.\@ Corti, A.\@ Pukhlikov and M.\@ Reid,
{\itshape Fano $3$-fold hypersurfaces},
Warwick, preprint, (1999).
\bibitem[Cu88]{Cu88}S.\@ Cutkosky,
{\itshape Elementary contractions of Gorenstein threefolds},
Math. Ann. 280 (1988), 521-525.
\bibitem[K88]{K88}Y.\@ Kawamata,
{\itshape Crepant blowing-up of $3$-dimensional canonical singularities
and its application to degenerations of surfaces},
Ann. Math. 127 (1988), 93-163.
\bibitem[K96]{K96}Y.\@ Kawamata,
{\itshape Divisorial contractions to $3$-dimensional terminal
quotient singularities},
Higher-dimensional complex varieties (Trento 1994), de Gruyter
(1996), 241-246.
\bibitem[KMM87]{KMM87}Y.\@ Kawamata, K.\@ Matsuda and K.\@ Matsuki,
{\itshape Introduction to the minimal model problem},
Adv. St. Pure Math. 10 (1987), 283-360.
\bibitem[$\mathrm{K}^+$92]{K+92}J.\@ Koll\'ar et al,
{\itshape Flips and abundance for algebraic threefolds},
Ast\'erisque 211 (1992).
\bibitem[M82]{M82}S.\@ Mori,
{\itshape Threefolds whose canonical bundles are not numerically effective},
Ann. Math. 116 (1982), 133-176.
\bibitem[M88]{M88}S.\@ Mori,
{\itshape Flip theorem and the existence of minimal models for $3$-folds},
J. Amer. Math. Soc. 1 (1988), 117-253.
\bibitem[R87]{R87}M.\@ Reid,
{\itshape Young person's guide to canonical singularities},
Proc. Symp. Pure Math. 46 (1987), 345-414.
\end{thebibliography}

\textsc{Department of Mathematical Sciences, University of Tokyo, Komaba,
Meguro, Tokyo 153-8914, Japan} \\
\texttt{kawakita@ms.u-tokyo.ac.jp}
\end{document}